\documentclass{amsart}
\usepackage[T1]{fontenc}
\usepackage[left=2cm, right=2cm, top=3cm, bottom=3cm]{geometry}

\usepackage{graphicx}

\usepackage{mathtools}
\usepackage{float}
\usepackage{booktabs}
\usepackage{siunitx} 
 \usepackage{multirow} 
\usepackage{subcaption}
\usepackage{comment}
\usepackage{thmtools}
\usepackage{booktabs}
\usepackage{siunitx}
\usepackage[table]{xcolor}

\DeclareMathOperator{\ESS}{ESS}

\usepackage{amsmath, amssymb}
\usepackage{algorithm}
\usepackage[noend]{algpseudocode}
\definecolor{lightgray}{rgb}{0.95,0.95,0.95}

\usepackage{nameref}
\usepackage{babel}
\usepackage{hyperref}

\DeclareMathOperator{\diff}{d}

\title{Data assimilation using a global Girsanov nudged particle filter}

\author{Maneesh Kumar Singh$^1$}
\author{Joshua Hope-Collins$^1$}
\author{Colin J. Cotter$^{1,*}$}
\author{Dan Crisan$^1$}
\address{$^1$ Department of Mathematics, Imperial College London, United Kingdom}
\address{$^*$ Corresponding author}
\date{}

\begin{document}

\begin{abstract}
    We present a particle filtering algorithm for stochastic models on infinite dimensional state space, making use of 
    Girsanov perturbations to nudge the ensemble of particles into regions of higher likelihood. We argue
    that the optimal control problem needs to couple control variables for all of the particles to maintain
    an ensemble with good effective sample size (ESS). We provide an optimisation formulation that separates the
    problem into three stages, separating the nonlinearity in the ESS term in the functional with the nonlinearity due to the forward problem, and allowing independent parallel computation for each particle when calculations
    are performed over control variable space. The particle filter is applied to the stochastic Kuramoto-Sivashinsky 
    equation, and compared with the temper-jitter particle filter approach. We observe that whilst the nudging filter is over spread compared to the temper-jitter filter, it responds to extreme events in the assimilated data more quickly and robustly.
\end{abstract}

\maketitle

\section{Introduction}

Data assimilation can be defined as the blending of observational data and computational models to obtain estimates of the past, present and future state of a physical system. It is a critical
tool in weather and climate which is now spreading across the engineering, physical and biological sciences \cite{reich2015probabilistic,nakamura2015inverse,law2015data,carrassi2018data}. In this article we focus on the filtering problem, which seeks 
the Bayesian estimate of the probability distribution for the model solution at time $t=t_0$, conditioned on noisy observations at times $0<t\leq t_0$ and a prior distribution at time $t=0$, and 
its solution by sequential methods that incrementally increase $t_0$. This distribution can be used to quantify uncertainty about the current state as well as the future state by forward propagation
using the computational model, to produce a forecast. The gold standard for approximate computational solution of the filtering problem is the particle filter, which evolves 
an empirical approximation of the filtering distribution using an ensemble of model states that are propagated in time using the computational model and then conditioned on the observed data 
by weighting and resampling. For the classical bootstrap particle filter, this ensemble needs to be large for useful estimates. A summary of analytical results for particle filters is given in \cite{crisan2002survey}.

When the model state space is large, as arising from the discretisation of a (possibly stochastic) partial differential equation (PDE), having a large number of particles  becomes impractical, and modifications and compromises are sought. In operational systems, the ensemble Kalman filter is a highly successful approach, based on approximating the distributions as multivariate normal and calculating a transformation of the ensemble of states instead of weighting and resampling \cite{evensen2009data}. These are combined with localisation
techniques to avoid finite size effects causing spurious long range correlations \cite{houtekamer1998data,hunt2007efficient}. However, there has also been a number of efforts to adapt particle filters to the smaller ensemble setting, such as equal weights particle filters \cite{ades2015equivalent,zhu2016implicit}, filters with Gaussian resampling \cite{potthast2019localized}, filters with ``tempering'' and ``jittering'' \cite{beskos2017stable,cotter2020particle,cotter2020data}, and filters with localisation \cite{poterjoy2016localized,crisan2025localized}. The goal of these efforts is to avoid assumptions about the distribution that underpin the ensemble Kalman filter.

In this article we focus on the incorporation of a ``nudging'' procedure into particle filters. Nudging is a long established technique in data assimilation \cite{kalnay2003atmospheric} that simply adds a relaxation term that ``steers'' the solution towards the observed data. More recently, there has been a programme of rigorous analysis to prove convergence of variants of these nudging schemes in the large time limit \cite{olson2003determining,azouani2014continuous,bessaih2015continuous}. In particle filters, the goal of nudging is to move the ensemble of states such that their likelihood is increased, without losing consistency of the filter. One approach is to modify the proposal \cite{chorin2010implicit,zhu2016implicit}. On the other hand, it was shown in \cite{akyildiz2020nudging,gonzalez2025nudging} that ensemble members can simply be nudged in directions of increasing likelihood which introduces bias but still gives asymptotic convergence for large ensembles provided that certain assumptions are met. In this article, we take a different approach based on a  Girsanov modification of a stochastic (partial) differential equation (S(P)DE). In this setting, the modification adds deterministic control variables to the S(P)DE
with a corresponding reweighting specified by Girsanov's Theorem. We are free to choose these variables provided that they do not depend on the forward filtration of the stochastic process. This means that they can only depend on the noise realisations for times in the past. In this work, we propose to use the control variables to nudge the ensemble towards regions of high Girsanov-modified likelihood, whilst ensuring that the ensemble members are roughly equally weighted. This can be viewed as a multiobjective optimisation problem. Solving this problem couples together the control variables for all of the ensemble members, which can be highly computationally challenging. In this article, we provide an approximate solution of this problem by breaking it down into three stages. The first and third stages solve uncoupled problems for each ensemble member, which can be solved in parallel. These problems depend on the forward solution of the S(P)DE, and their practical solution is enabled by automated adjoint computation. The second stage connects all of the ensemble members together, but only requires the solution of a constrained optimisation problem in $\mathbb{R}^N$, where $N$ is the ensemble size, and without requiring forward/adjoint model solutions. Hence, the algorithm is feasible in the high performance computing setting. 

The rest of this paper is organised as follows. In Section \ref{sec:pf} we establish some notation and preliminaries before describing our nudging framework and the three stage solution approach. In Section
\ref{sec:numex} we provide some numerical experiments illustrating our approach, and in Section \ref{sec:summary} we provide a summary and outlook.

\section{Girsanov nudged particle filtering}
\label{sec:pf}
In this section, the details of our particle filtering approach will be presented.

\subsection{Notation and preliminaries}

We adopt the stochastic nonlinear filtering framework for data assimilation. Let $(X, Y)$ be processes on $(\Omega, F, \mathbb{P})$, where $X$ is the (hidden) signal and $Y$ is a discrete time observation process 
\begin{equation}\label{obseq}
    Y_{t_k} := h(X_{t_k})) + V_{t_k},
\end{equation}
where $h: \Omega \to \mathbb{R}^M$ is 
the observation operator, $V_t$ is a discrete time
noise process, and $(t_1,t_2,t_3,\ldots)$ is a sequence
of observation times. The probability density function $L(Y_{t_k}|X_{t_k})$ for
$Y_{t_k}$ conditional on the value of $X_{t_k}$ is referred to as the likelihood.
We will work with the negative log likelihood $\Phi(X_{t_k}, Y_{t_k})=-\log(Y_{t_k}|X_{t_k})$.

Here, $X$ solves a given SPDE, and the goal is to estimate
the filtering distribution $\pi_{t_k}(X_{t_k}|Y_{t_1},Y_{t_2},\ldots,Y_{t_k})$ for
the signal at $t=t_k$ given all observations $Y_{t_j}$
for $j=1,\ldots,k$. Bayes' rule states that
\begin{equation}
\frac{\diff\pi_{t_k}(X_k|Y_{t_1},Y_{t_2},\ldots,Y_{t_{k-1}},Y_{t_k})}{\diff \pi_{t_k}(X_k|Y_{t_1},Y_{t_2},\ldots,Y_{t_{k-1}})} = L(Y_{t_k}|X_{t_k}),
\end{equation}
where $d\pi_A/d\pi_B$ is the Radon-Nikodym derivative of $\pi_A$ with respect to $\pi_B$.

We approximate $\pi_t$ using \textit{particle filters}, which are 
sequential Monte Carlo methods representing the conditional law of $X_t$ via a weighted ensemble of model states, or particles, $\{(w_i,X_i)\}_{i=1}^{N_p}$, with
\[
\pi_{t} \approx \sum_{i=1}^{N_p} w_i\delta(X-X_i(t)).
\]
In the most basic bootstrap particle filter, particles are propagated forward in time 
from $t=t_{k-1}$ to $t=t_k$ using independent realisations of the stochastic process modelling the state
for each particle. Then, each particle $X_{i,t}$ is assigned an updated likelihood weight \emph{via}
\[
w_i \mapsto \hat{w}_i = w_i L(X_{i,t_k}|Y_{t_k}),
\]
which are then normalised according to $w_i=\hat{w}_i/\sum_{k=1}^{N_p}\hat{w}_k$. 

To assess weight degeneracy, we compute the \textit{effective sample size},
\begin{equation}\label{ess}
\mathrm{ESS}({\mathbf{w}}) \coloneqq \left( \sum_{i=1}^{N_p} {w}_i^2 \right)^{-1}.
\end{equation}
Here, $\mathrm{ESS} \approx N_p$ implies uniform weights, while $\mathrm{ESS} \ll N_p$ indicates a loss of ensemble spread which will lead to filter divergence. Resampling is triggered when the normalised ESS, $N_p^* = \mathrm{ESS}/N_p$, falls below a threshold, resulting
in a new set of equally weighted particles which are sampled from the original set with replacement in such 
a way that the empirical distributions before and after resampling are asymptotically consistent as $N_p\to \infty$.
This has the effect of replicating higher weight particles and discarding lower weight particles (on average).
In this article we used the systematic resampling algorithm of \cite{gordon1993novel}.

\subsection{The Girsanov nudging framework}

The effective sample size (ESS) reflects weight degeneracy, and we frequently observe it dropping in high dimensions (such as discretisations of SPDEs) due to particle sparsity in observation space. This leads to poor posterior approximation. 

In this framework, we introduce a time-dependent \textit{control variable} to guide particles toward regions suggested by observations. This \textit{nudging} preserves consistency, meaning the particle ensemble remains a valid set of samples from the prior distribution after appropriate weight modification. 

We justify this approach through  a stochastic differential equation (SDE),
\begin{equation}
\label{eq:sde}
\diff X = f(X)\diff t + \sigma(X)\diff W, \quad X(0)=X_0,
\end{equation}
where $X \in \mathbb{R}^{N_1}$, $f:\mathbb{R}^{N_1}\to \mathbb{R}^{N_1}$, $\sigma: \mathbb{R}^{N_1}\to \mathbb{R}^{N_1 \times N_2}$, and $W(t)$ is a $N_2$-dimensional Brownian motion. If we instead solve a modified SDE,
\begin{equation}
\label{eq:gsde}
\diff \hat{X} = f(\hat{X})\diff t + \sigma(\hat{X})\left(\lambda(t)\diff t + \diff W\right), \quad \hat{X}(0) = X_0,
\end{equation}
where $\lambda(t) \in \mathbb{R}^{N_2}$, the joint probability measure for $\hat{X}(t)$ is absolutely continuous with respect to that of $X(t)$. The Radon-Nikodym derivative (or Girsanov factor) between these measures is given by
\begin{equation}
\label{eq:girsanov factor}
G = \exp\left(-\int_0^T \frac{1}{2}|\lambda(t)|^2 \diff t - \int_0^T \lambda(t)\cdot\diff W\right),
\end{equation}
subject to appropriate regularity conditions on $f(X)$ and $\sigma(X)$, and subject to the condition that
$\lambda(t)$ does not depend on $W(s)-W(t)$ for $t>s$. 
Thus, we can choose $\lambda$ to steer $X$ towards regions of higher likelihood function, but we must compensate by multiplying this weight by $G$. In other words, we pay a cost for nonzero $\lambda(t)$.

In the discrete time setting, we need to approximate the time integrals in \eqref{eq:girsanov factor} by
quadrature sums, e.g.,
\begin{equation}
\label{eq:girsanov DG approx}
G \approx G_{\Delta t} = \exp\left(
-\sum_{n=1}^{N_s}\left(\frac{1}{2} \lambda^n)^2 \Delta t 
- \lambda^n\Delta W^n\right)\right),
\end{equation}
where $\Delta W^n$ are suitable independent Brownian increments and $\Delta t$ is the time stepsize.

In the context of particle filtering, this means that we can solve \eqref{eq:gsde} instead of \eqref{eq:sde}
to propagate each particle from $t_{k-1}$ to $t_{k}$, but now the weight must be updated using the Girsanov-adjusted
formula,
\[
w_i \mapsto \hat{w}_i = w_i \exp(-\hat{\Phi}_i),
\]
where
\begin{equation}
\hat{\Phi}_i = \Phi(X_{i,t_k},Y_{t_k})+\int_{t_{k-1}}^{t_k} \frac{1}{2}\|\lambda_i(t)\|^2 \diff t - \int_{t_{k-1}}^{t_k} \lambda_i(t)\cdot\diff W, \quad i=1,\ldots, N_p.
\end{equation}
Here we can select different controls $\lambda_i$ for each particle $X_i$.

This leads naturally to the idea of maximising $\hat{\Phi}_i$ with respect to $\lambda_i$. However, it is essential to ensure that Girsanov’s theorem remains valid. This requires that the control $\lambda(t)$ depends only on the past values of the noise process, i.e., on $W(s)$ for $s < t$. To satisfy this adaptivity constraint, we adopt a sequential optimisation strategy. The idea is to construct $\lambda(t)$ incrementally as we uncover new values of $W(t)$. Upon discretising time, and defining $\lambda_{i}(t_{k-1} + n\Delta t) = \sum_{j=1}^i \Delta \lambda_{j,n}, i=1,2,\ldots,N_p$,
$n=1,2,\ldots,N_s$ such that $t_k=t_k+N_s\Delta t$, 
we proceed as follows. First, initialise $\Delta\lambda_{i,n} = 0$ and $\Delta W_i^n = 0$ for $n = 1, 2,\ldots, N_s$,
and $i=1,\ldots,N_p$. 
Then we sample $\Delta W_{i,1}$ from the appropriate distribution, and adjust $\Delta\lambda_{i,1}$ to decrease
$\hat{\Phi}_i$, when computed by solving the forward model for $X_i$, with all other parameters fixed,
for $i=1,2,\ldots,N_p$. Then we sample $\Delta W_{i,2}$ and adjust $\Delta\lambda_{i,2}$ and so on, until we reach
$n=N_s$.

\subsection{Three stage optimisation process}
In previous work \cite{cotter2020data, cjc_da_mks_SCH} we have investigated the application of this idea to particle filtering for SPDEs. The value of $\lambda_i$ was chosen to minimise
$\hat{\Phi}_i$, thus maximising the Girsanov-adjusted weight. This amounts to resolving a trade-off between
maximising the likelihood, and minimising the Girsanov ``cost''. \cite{cotter2020data} made an initial investigation where $\lambda_i$ was only nonzero for the final stage of 
the final timestep before the assimilation time $t_k$. This restriction means that $\lambda_i$ can be obtained as the solution of a least squares problem. However, this provides limited opportunity to steer the solution between assimilation times, only providing an impulsive kick right at the end. \cite{cjc_da_mks_SCH} implemented a more flexible optimisation framework where $\hat{\Phi}_{i}$ is optimised for more general time dependencies in $\lambda_i$
(stepwise constant values in this case), along with a parallel implementation where the optimisation problems are
solved independently in parallel. These optimisation problems are solved using a gradient method, with the gradient being computed \emph{via} the adjoint equations. In that work, we established that this approach can produce consistent filtering solutions, but it required combination with additional tempering and jittering steps as introduced in \cite{beskos2017stable,cotter2020particle} to obtain stable results. This is because independent
optimisation of $\hat{\Phi}_i$ for each particle $X_i$ does not address the problem of filter divergence. The frequent
failure of the bootstrap filter in high dimensions is because a small number of particles will acquire almost all
of the normalised weight, because they have larger likelihoods than the others and the likelihood function is typically exponentially decaying. This is observable in low ESS values, and results in almost all particles having
the same state after resampling. If we individually optimise each $\hat{\Phi}_i$, although all particles are moved to
regions of higher likelihood, the particles that would have otherwise been close can be moved even closer, and the 
ESS remains low.

To try to keep ESS as high as possible, we propose to
solve a global optimisation problem over all particles, \emph{i.e.} over the full set $(\lambda_1,\lambda_2,\ldots,\lambda_{N_p})$ of control variables. This is closely related to the idea behind
the equal weights particle filter \cite{zhu2016implicit}. Here, we consider maximising the functional
\begin{equation}
ESS = \frac{\left(\sum_{i=1}^{N_p} \exp(-\hat{\Phi}_i)\right)^2}{\sum_{i=1}^{N_p} \exp(-2\hat{\Phi}_i)},
\end{equation}
which reduces to the usual ESS formula upon setting $\hat{w}_i=\exp(-\hat{\Phi}_i)$ and applying normalisation.
(Note that care must be taken here to avoid issues with numerical underflow.) Due to the normalisation factor,
this optimisation problem is ill-posed, since scaling all the weights by a constant will not change the ESS.
To deal with this we add a regularisation term, minimising 
\begin{equation}
\label{eq:F}
\mathbb{F}_{\ESS} = \sigma \sum_{i=1}^{N_p} \hat{\Phi}_i - \dfrac{\left(\sum_{i=1}^{N_p}\exp(-\hat{\Phi}_i)\right)^2}{\sum_{i=1}^{N_p} \exp(-2\hat{\Phi}_i)},
\end{equation}
where $\sigma>0$ is a chosen penalty coefficient.
This modification will tend to move particles into regions of Girsanov adjusted weight if it is possible without unduly reducing ESS.

In our experiments, optimisation algorithms for this problem converged very slowly, probably due to the composition of the nonlinearity of the forward model producing $\hat{\Phi}_i$ with the nonlinearity in the ESS. This problem
also requires an additional layer of parallel communication between particles which are otherwise running independently. To circumvent this, we designed the following three stage algorithm to find approximate solutions
of this optimisation problem. In the first stage, we find $\lambda_i$ that minimises $\hat{\Phi}_i$ separately
for each particle $X_i$. This can be done for each particle independently in parallel, and determines the possible ranges of values $\hat{\Phi}_{i,\min}$ (from the optimising $\lambda_i=\lambda_i^*$) and $\hat{\Phi}_{i,\max}$ (from $\lambda_i=0$) for $\hat{\Phi}_i$ in the global optimisation. In the second stage, we minimise $\mathbb{F}_{\ESS}$ in \eqref{eq:F} over $(\hat{\Phi}_1,\hat{\Phi}_2,\ldots,\hat{\Phi}_{N_p})$ instead of over $(\lambda_1,\lambda_2,\ldots,\lambda_{N_p})$, subject to the
constraints 
\[
\hat{\Phi}_{i,\min} \leq \hat{\Phi}_{i} \leq \hat{\Phi}_{i,\max}, \quad i=1,2,\ldots,N_p.
\]
This is a problem in state space dimension $N_p$, and does not involve running the forward model.
After computing the optimal values $(\hat{\Phi}^*_1,\hat{\Phi}^*_2,\ldots,\hat{\Phi}^*_{N_p})$ in stage two, in stage three we
need to find the value of $\lambda_i$ that produces $\hat{\Phi}_i=\hat{\Phi}^*_i$ for $i=1,2,\ldots,N_p$.
There are multiple solutions to this problem, but we set $\lambda_i=s_i\lambda_i^*$ for $s_i\in \mathbb{R}$,
and seek the value of $s_i$ that achieves the optimal $\hat{\Phi}^*_i$, which is reachable provided that 
$\hat{\Phi}_i$ is continuous when viewed as a function of $s_i$. This is a one-dimensional nonlinear equation that can also be solved independently
in parallel for each particle $X_i$.

To state the algorithm more precisely, we assume that the SPDE is discretised such that there is 
just one (multidimensional) Brownian increment sampled per timestep $n$, $\Delta W_{i,n}$, which
we also label with the particle index $i$ to indicate that dependency. Then, we define
$F(X_{i,t_{k-1}}, \Delta {W}_i, \Delta \lambda_i)$ as the solution obtained at $t=t_k$ from the 
numerical discretisation of the SPDE using $\Delta W_i=(\Delta W_{i,1},\Delta W_{i,2},\ldots,\Delta W_{i,N_s})$
and $\Delta \lambda_i=(\Delta \lambda_{i,1},\Delta \lambda_{i,2},\ldots,\Delta \lambda_{i,N_s})$.
We also indicate dependency of the Girsanov-adjusted negative log likelihood values via
\[
\hat{\Phi}(X_{i,t_k}, Y_{t_k}, \Delta W_i, \Delta \lambda_i).
\]
The values of $\Delta W_i,\Delta \lambda_i$, $i=1,2,\ldots,N_p$, are then updated according to
Algorithm \ref{al:tsn}.

\begin{algorithm}
\caption{\label{al:tsn} Three Stage Nudging}
\begin{algorithmic}
\State Set $\Delta W_i=0$, $\Delta \lambda_i=0$.
\For{$n=1$ to $N_s$}
 \State \textbf{Stage 1: Discover minimum bounds for $\hat{\Phi}$}
 \For{$i=1$ to $N_p$ (in parallel)}
 \State Sample $\Delta W_{i,n}$
 \State Set $\hat{\Phi}_{i,\max}=\hat{\Phi}(F(X_{i,t_{k-1}}, \Delta W_i, \Delta \lambda_i), Y_{t_k},
    \Delta W_i, \Delta \lambda_i)$
 \State
    Keeping all other parameters fixed, set 
    \[
    \Delta \lambda_{i,n}^*=\arg\max_{\Delta\lambda_{i,n}}\hat{\Phi}(F(X_{i,t_{k-1}}, \Delta W_i, \Delta \lambda_i), Y_{t_k},
    \Delta W_i, \Delta \lambda_i).
    \] 
 \State Set $\hat{\Phi}_{i,\min}$ as the obtained corresponding value of $\hat{\Phi}$.
 \EndFor
\State \textbf{Stage 2: Maximise ESS}
 \State Set
 \[
 (\hat{\Phi}^*_1,\hat{\Phi}^*_2,\ldots,\hat{\Phi}^*_{N_p})
 = \arg\min_{(\hat{\Phi}_1,\hat{\Phi}_2,\ldots,\hat{\Phi}_{N_p})} \left(-\text{ESS}(\hat{\Phi}_1,\hat{\Phi}_2,\ldots,\hat{\Phi}_{N_p}) + \sigma \sum_i \hat{\Phi}_i \right)
 \]
subject to 
\[\hat{\Phi}_i^{\min} \leq \hat{\Phi}_i^{*} \leq \hat{\Phi}_i^{\max}, \quad i=1,2,\ldots,N_p.
\]
\State \textbf{Stage 3: Recover Scaling Factor}
\For{$i=1$ to $N_p$ (in parallel)}
\State Find $s$ such that $\Phi_i(X_{i,t_{k-1}},Y_{t_k},\Delta W_i, \widehat{\Delta \lambda}_i) = \Phi_i^*$
where $\widehat{\Delta \lambda}_{i,n} = s\Delta \lambda_{i,n}$ and $\widehat{\Delta \lambda}_{i,m} = \Delta \lambda_{i,m}$ for $m\ne n$.
\State Scale $\Delta \lambda_{i,n}$ by $s$.
\EndFor
\EndFor
\end{algorithmic}
\end{algorithm}

\section{Numerical experiments}
\label{sec:numex}

In this section, we demonstrate the Girsanov nudging approach with some numerical experiments. These were 
implemented using our particle filter library \cite{nudging} which is built around Firedrake \cite{FiredrakeUserManual}, an automated system for solving partial differential equations using the finite 
element method. Our library has a modular design which allows for any SPDE discretisation that is expressible
in Firedrake; new discretisations are added by subclassing a model class which hands over noise realisation sampling to the particle filter so that different filter algorithms can be used seamlessly. The library makes use of
Firedrake's automated adjoint capability \cite{farrell2013automated,mitusch2019dolfin}, allowing gradients of $\hat{\Phi}$ with respect to $\Delta \lambda$ to be computed scalably and nonintrusively in the optimisation procedure. The library also makes use of message passing interface (MPI) subcommunicators so that the algorithm
can be parallelised over particles as well as composed with standard spatial domain decomposition parallelism applied
to the forward model. In this article, we only demonstrate the scheme with  SDE and  one-dimensional SPDEs, where only parallelism 
over particles was used, but the system is designed to be fully capable for large scale forward models, which we will consider in the future. 

In our numerical experiments, we compared our particle filter with (a) the classical bootstrap particle filter,
and (b) the temper-jitter particle filter of \cite{beskos2017stable,cotter2020data,cotter2020particle}. This
latter filter transforms from the prior distribution $\pi_{t_k^-}$ to the posterior distribution $\pi_{t_k^+}$
in a sequence of steps via intermediate distributions, $\pi_{\theta_j}$, $j=0,\ldots,N_\theta$, with $0=\theta_0=\theta_1<\ldots<\theta_{N_\theta}:=1$, and
\[
\frac{\diff{\pi_{\theta_j}}}{\diff\pi_{\theta_{j-1}}} = \left(L(X_{k},Y_k)\right)^{\Delta \theta_j}, \quad \pi_{\theta_0}=\pi_{t_k^-}, \quad \Delta \theta_j = \theta_j-\theta_{j-1}.
\]
At each step, the temper-jitter filter selects $\Delta \theta_j$ adaptively so that the ESS remains above 
some threshold ($ESS > 0.8N_p$ in our experiments). The particles are weighted by the ``tempered likelihood'' $\left(L(X_{k},Y_k)\right)^{\Delta \theta_j}$, and then resampled to obtain an equal weighted ensemble again.
Then, ``jittering'' is applied, meaning that some number of Monte Carlo Markov Chain (MCMC) iterations are applied
to the noise increments used to advance the particle state from $t_k$ to $t_{k-1}$, given the initial condition
$X_{i,k-1}$, to obtain alternative samples from the distribution $\pi_{\theta_j}$. This has the effect of avoiding
duplication of particles and improving the approximation, by selecting alternative noise increments that are statistically consistent with the target distribution. The jittering steps are applied independently in parallel for each particle and the only coupling between particles occurs in the resampling step. In our experiments
we used the Preconditioned Crank Nicholson MCMC proposal of \cite{cotter2013mcmc} in our jittering steps,
in the form of their Equation (4.6) for a chosen step parameter $\delta>0$.
This process is then repeated to transform from $\pi_{\theta_k}$ to $\pi_{\theta_{k+1}}$ and so on. This process requires to continue to store the initial
conditions $X_{i,k-1}$ and noise increments $dW_{i}$ during the resampling process.

For the experiments with the Girsanov nudging particle filter, we found that applying a small number of jittering
steps (without tempering) after resampling improved the quality of the results in the case of the linear SDE
example. Hence, we have done this in all of our results. In these experiments, we were mainly focused on the 
behaviour of the filter in terms of accuracy and stability, and we did not attempt to tune termination criteria
for the optimisation solvers, leaving this for future work with more advanced versions of the filter as
discussed in the summary and outlook section.


\subsection{Linear SDE}
First, as a brief verification of consistency of the filter, we 
consider the one-dimensional linear SDE
\begin{equation}
\label{eq:lsde}
\mbox{d}  x = -A x  \mbox{d}  t + D\mbox{d}W,
\end{equation}
for $x\in \mathbb{R}$, with $A,D\in \mathbb{R}$ positive constants. The corresponding
Girsanov-perturbed SDE is
\begin{equation}
\label{eq:glsde}
\mbox{d}x = -A x  \mbox{d}  t + D(\mbox{d}W + \lambda \mbox{d}  t ).
\end{equation}
We discretise \eqref{eq:glsde} using the midpoint scheme,
\[
(1+A \Delta t/2)x^{n+1} = (1 - A \Delta t/2)x^{n} + D (\Delta W^{n+1} + \Delta t \lambda^{n+1}),
\]
where $\Delta W^{n+1}\sim N(0, \Delta t)$, and
where the scheme for \eqref{eq:lsde} is recovered from the choice $\lambda^{n+1}=0$.

In our experiment, we initialise the ensemble as samples from the $N(0, D^2/(2A))$ distribution, 
selecting $A=D=1$. The observation operator is $h(x)=x$, and we make a noisy observation
$h(x)+\epsilon=-0.055634$ with $\epsilon\sim N(0, 0.01)$ at time $t=1$, with time stepsize
$\Delta t = 1/10$. 

The Girsanov nudging particle filter was used with 5 jittering steps with MCMC step parameter $\delta=0.05$ 
after resampling. The Stage 1 optimisation was solved using L-BFGS (Limited memory version of the Broyden-Fletcher-Goldfarb-Shanno algorithm) \cite{nocedal2006numerical}, using the SciPy \cite{scipy} implementation as wrapped in the pyadjoint 
library \cite{mitusch2019dolfin} with parallelism over the particles. The Stage 2 optimisation was solved using the L-BFGS-B (bounded version of L-BFGS) provided by SciPy, and the Stage 3 optimisation was solved using Brent's method \cite{brent2013algorithms} also provided by SciPy, also with parallelism over the particles.

We compared the Girsanov nudging filter with the bootstrap filter and the temper-jitter filter, the latter of 
which used MCMC step parameter $\delta=0.15$ and 5 jittering steps after each tempered resampling step.

We compare the mean and the variance of the posterior distribution with the exact values for the undiscretised
SDE (hence the filter results will contain timestepping errors). Results are presented in Tables
\ref{tab:mean_performance_compact} and \ref{tab:variance_performance_compact}; we conclude that the nudging particle
filter produces estimates that are of similar accuracy to the temper-jitter filter. Remaining errors appear to
be due to the time discretisation of the forward model.

\begin{table}[htb!]
    \centering
    \footnotesize
    \setlength{\tabcolsep}{4pt} 
    \renewcommand{\arraystretch}{0.9} 
    \sisetup{
        table-align-text-post = false,
        detect-weight = true,
        detect-family = true
    }
    \begin{tabular}{
        l
        c
        S[table-format=2.0]
        S[table-format=-1.6]
        S[table-format=-1.6]
        S[table-format=1.6]
    }
        \toprule
        Filter & Ensemble size & {ESS (\%)} & {Exact mean} & {Ensemble mean} & {Estimator error} \\
        \midrule
        \rowcolor{lightgray} Bootstrap & 90  & 21 & -0.054543 & -0.071101 & 0.016558 \\
                                & 150 & 19 & -0.054543 & -0.050835 & 0.003709 \\
        \rowcolor{lightgray}    & 300 & 17 & -0.054543 & -0.042207 & 0.012336 \\
        \midrule
        \rowcolor{lightgray} Temper-Jitter & 90  & {}   & -0.054543 & -0.051902 & 0.002641 \\
                                     & 150 & {}   & -0.054543 & -0.066020 & 0.011477 \\
        \rowcolor{lightgray}         & 300 & {}   & -0.054543 & -0.064980 & 0.010437 \\
        \midrule
        \rowcolor{lightgray} Nudge(+Jitter) & 90  & 55 & -0.054543 & -0.057870 & 0.003327 \\
                                        & 150 & 51 & -0.054543 & -0.061563 & 0.007019 \\
        \rowcolor{lightgray}           & 300 & 51 & -0.054543 & -0.057404 & 0.002861 \\
        \bottomrule
    \end{tabular}
        \caption{Filter estimates of the posterior mean for the linear SDE experiment. The estimates are 
        displayed for the bootstrap filter, the temper-jitter filter, Girsanov nudging filter with jittering, alongside the exact value, for various ensemble sizes. The ratio of ESS to $N_p$ (ESS $(\%)$) is also shown for bootstrap and the Girsanov nudging filter, showing that the nudging is raising the ESS relative to the bootstrap filter. ESS is not shown for the temper-jitter filter since the tempering is adapted to ensure ESS $(\%) \approx 80$ at each tempering step.}
    \label{tab:mean_performance_compact}
\end{table}

\begin{table}[htb!]
    \centering
    \footnotesize
    \setlength{\tabcolsep}{4pt} 
    \renewcommand{\arraystretch}{0.9} 
    \sisetup{
        table-align-text-post = false,
        table-format = 1.6,
        detect-weight = true,
        detect-family = true
    }
    \begin{tabular}{
        l
        c
        S[table-format=1.6]
        S[table-format=1.6]
        S[table-format=1.6]
    }
        \toprule
        Filter & Ensemble size & {Exact variance} & {Ensemble variance} & {Estimator error} \\
        \midrule
        \rowcolor{lightgray} Bootstrap & 90  & 0.009804 & 0.008484 & 0.001320 \\
                                & 150 & 0.009804 & 0.007480 & 0.002324 \\
        \rowcolor{lightgray}    & 300 & 0.009804 & 0.008547 & 0.001257 \\
        \midrule
        \rowcolor{lightgray} Temper-Jitter & 90  & 0.009804 & 0.013017 & 0.003213 \\
                                     & 150 & 0.009804 & 0.010127 & 0.000323 \\
        \rowcolor{lightgray}         & 300 & 0.009804 & 0.009823 & 0.000019 \\
        \midrule
        \rowcolor{lightgray} Nudge(+Jitter) & 90  & 0.009804 & 0.010390 & 0.000586 \\
                                        & 150 & 0.009804 & 0.009780 & 0.000024 \\
        \rowcolor{lightgray}           & 300 & 0.009804 & 0.009361 & 0.000442 \\
        \bottomrule
    \end{tabular}
        \caption{Filter estimates of the posterior variance for the linear SDE experiment, with the same formatting as Table \ref{tab:mean_performance_compact}.}
    \label{tab:variance_performance_compact}
\end{table}

\subsection{Stochastic Kuramoto-Sivashinsky equation}

To investigate the stability and accuracy of the Girsanov nudging approach, we performed data assimilation
experiments using the stochastic Kuramoto-Sivashinsky (SKS) equation with additive noise \cite{duan2001stochastic},
\begin{equation}\label{skse1}
\diff u + \left(\alpha u_{xxxx} + \beta u_{xx} + \gamma uu_x\right)\diff t = c \diff W,
\end{equation}
solved on the interval $[0,L]$ with periodic boundary conditions $u(L,t)=u(0,t)$,
where $\alpha$, $\beta$, $\gamma$ and $c$ are constants, and where $W$ is a spacetime white noise
(cylindrical in space and Ito in time). In the absence of noise, this equation has a global attractor
and exhibits chaotic behaviour for suitable parameter choices \cite{papageorgiou1991route}. The stochastic
extension has random attractors with stable long time behaviour \cite{yang2006random}. This makes  the SKS model  very useful as a tool for rapidly benchmarking particle filters.

We discretised \eqref{skse1} using the degree 2 continuous Lagrangian finite element space $V_h$ on a regular mesh
with $N_v$ vertices and mesh width $h=L/N_v$,
using a $\mathcal{C}^0$ interior penalty (CIP) treatment of the fourth order term \cite{brenner2005c}.
For the time discretisation, we consider the mid-point scheme on a uniform time mesh.
The spatially discrete continuous time scheme is given by
\begin{equation}
\left(\diff u_h,v_h\right) - \left(\beta u_{h,x}\diff t, v_{h,x}\right) + a\left(\alpha u_{h}\diff t, v_h\right) 
-\left(\dfrac{\gamma}{2} \left(u_{h}\right)^{2}\diff t, v_{h,x}\right)=c\diff W_h[v_h], \quad
\forall v_h \in V_h,
\end{equation}
where
\begin{align*}
  a(u,v) &= (u_{xx},v_{xx}) + \langle \{u_{xx}\}, [v_x]\rangle+ \langle \{v_{xx}\}, [u_x]\rangle + \dfrac{\eta}{h}\langle [u_x], [v_x]\rangle, \\
  (u,v) & = \int_0^L uv\diff x, \\
  \langle u,v \rangle & = \sum_{i=1}^{N_v}u(z_i)v(z_i),\\
  \{u\}|_{z_i} & = (u(z_i^+)+u(z_i^-))/2, \\
  [u]|_{z_i} & = u(z_i^+)-u(z_i^-),
  \end{align*}
where  $^+$ and $^-$ indicates right and left limiting values 
  respectively, $\eta$ is the (mesh independent) interior penalty parameter, and
\begin{equation}
\label{eq:dW_h}
\diff W_h[v] = \frac{1}{h^{1/2}}\sum_i^{N_v} \left(\int_{ih}^{(i+1)h} v \diff x\right)\diff W_i,
\end{equation}
with $\{W_i\}_{i=1}^{N_v}$ a set of iid Brownian motions, so that
\begin{align}
\mathbb{E}\left(\left(\int_{t_0}^{t_1}\diff W_h[v]\right)
\left(\int_{t_0}^{t_1}\diff W_h[w]\right)\right)
&= \int_{t_0}^{t_1}\frac{1}{h}\sum_i^{N_v} 
\left(\int_{ih}^{(i+1)h} v \diff x\right)
\left(\int_{ih}^{(i+1)h} w \diff x\right)\diff t, \\
& \approx \int_{t_0}^{t_1}\int_0^L vw \diff x \diff t. \nonumber
\end{align}
We apply the Girsanov perturbation after discretising in space but before discretising in time, replacing
$\diff W_i\mapsto \diff W_i + \lambda_i(t)\diff t$, $i=1,2,\ldots,N_v$ in \eqref{eq:dW_h}. After applying
the implicit midpoint rule, 
the fully discretised Girsanov perturbed system is
\begin{align}
\nonumber
\left({u_h^{n+1}-u_h^n},v_h\right) -  \left(\alpha\Delta tu_{h,x}^{n+1/2}, v_{h,x}\right) + a\left(\Delta t \beta u_{h}^{n+1/2}, v_h\right) & \\
\quad -\left(\dfrac{\gamma \Delta t}{2} \left(u_{h}^{n+1/2}\right)^{2}, v_{h,x}\right)&=c\widehat{\Delta W_h}^n[v_h],
\quad \forall v_h\in V_h,
\end{align}
where
\begin{equation*}
\widehat{\Delta W_h}^n[v] = \frac{1}{h^{1/2}}\sum_i^{N_v} \left(\int_{ih}^{(i+1)h} v \diff x\right)
\left(\Delta W_i^n + \lambda_i^n \Delta t\right),
\end{equation*}
where $\Delta W_i^n$ are iid $N(0,\Delta t)$ random variables. Here we see that $\lambda$ takes the role 
of a piecewise constant function in space. The unperturbed discretised SPDE is recovered for
$\lambda_i^n=0$.

We employed a 90-particle ensemble in our numerical particle filtering studies applied to the SKS model. This is motivated by common ensemble sizes of data assimilation and ensemble uncertainty quantification systems for operational weather forecasting, where the high computational cost of the forward models limits ensemble sizes; our goal is to obtain stable particle filter configurations in the low ensemble size setting.

In the numerical setup of the SKS model, we chose the length of the spatial domain $L=4$ with 100 cells (so that there are 200 degrees of freedom, since the finite element space is piecewise quadratic), and we take the diffusion coefficient $\alpha = 0.03$, the anti-diffusion term $\beta = 1.1$, the advection parameter $\gamma = 1$, and the noise parameter $c = 2.5$. We used $\eta=5$ for the interior penalty parameter.

\begin{figure}[htb!]
        \centering
        \includegraphics[width=0.6\linewidth]{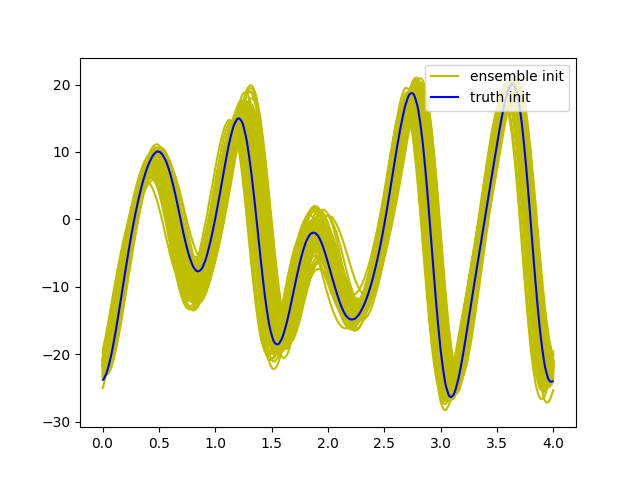} \\
        \includegraphics[width=0.6\linewidth]{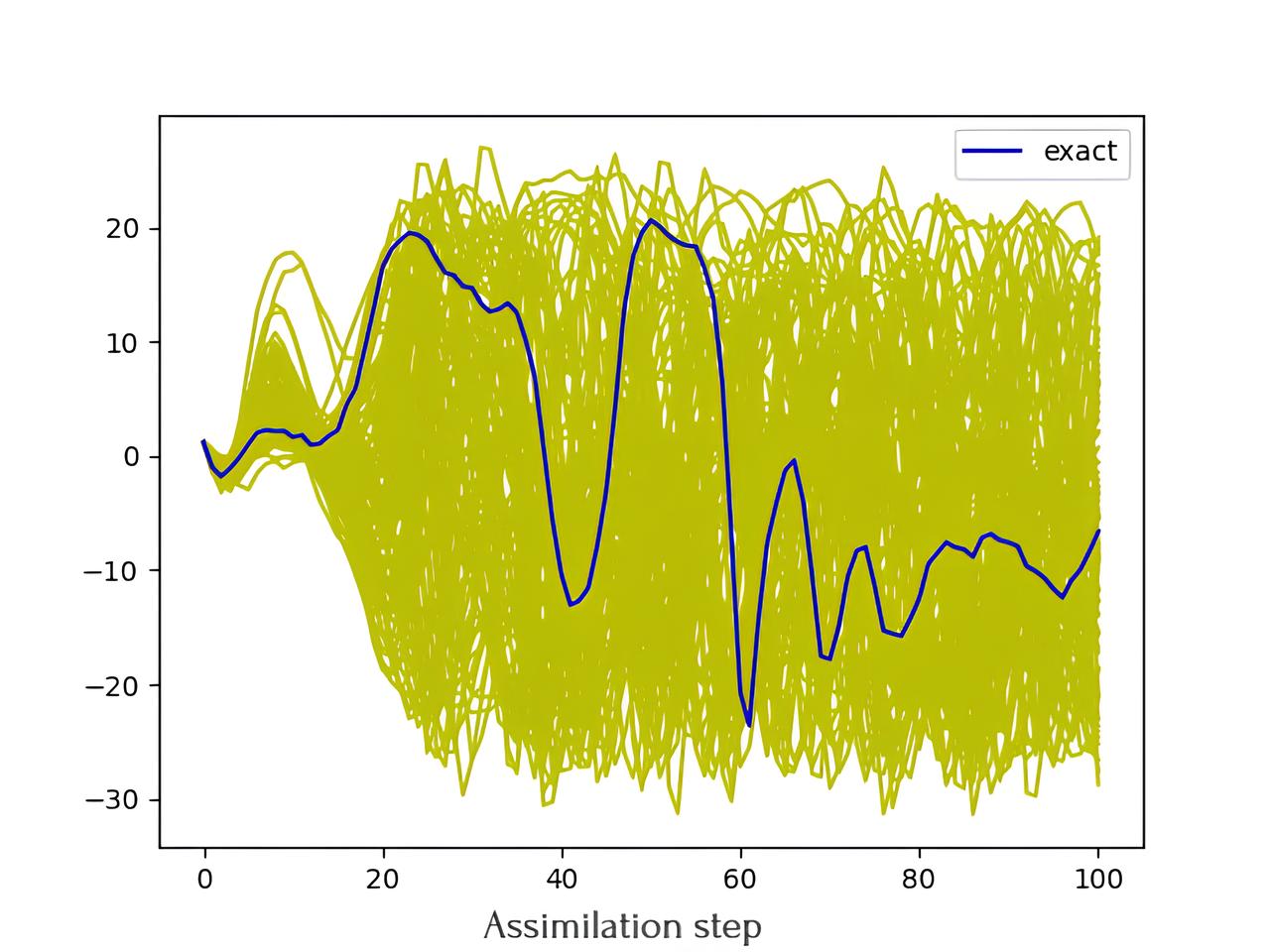}
        \caption{Top: Initial distribution of particles (yellow colour) for the stochastic Kuramoto-Sivashinsky system, with the ``true'' value used to generate observations shown in blue. Bottom: a visualisation
        of how the ensemble of particles spreads when the observations are not incorporated, showing the time 
        evolution of the solution at one observation point. The true value (indicated as "exact") is also shown.
        The time scale is shown in terms of assimilation steps to allow comparison with other Figures.}
        \label{fig:init}
\end{figure}

For all the numerical experiments, the initialisation of particles and reference trajectory (``truth'') is constructed in the following way.
First, we run the above numerical simulation of the SKS model  over 200 steps with the following initial condition,
\[
u_{in}(x) = \left(\dfrac{0.4}{e^{(x-403./15.)} + e^{(-x+403./15.)}} + \dfrac{1}{e^{(x-203./15.)}+e^{(-x+203./15.)}}\right).
\]
Then, the solution at the final time of the above simulation is considered as the initial condition $u_{0}$ for the truth and all the particles. Fig \ref{fig:init} displays the initialisation of truth and particles.  We observe that the particle distribution is somewhat spread around the truth. Furthermore, we have demonstrated  evolution of particles,  without data assimilation. As expected, a significant spread is observable when comparing these particles to the truth.

We took measurements at $ 10$ equispaced points in the
interval $[0, 4]$. All observation processes are perturbations of the ``true'' trajectory with iid measurement errors of distribution $\mathcal{N}(0, 2.5)$. The observations are taken every 5 timesteps.

\begin{figure}[htb!]
        \centering
        \includegraphics[width=0.5\linewidth]{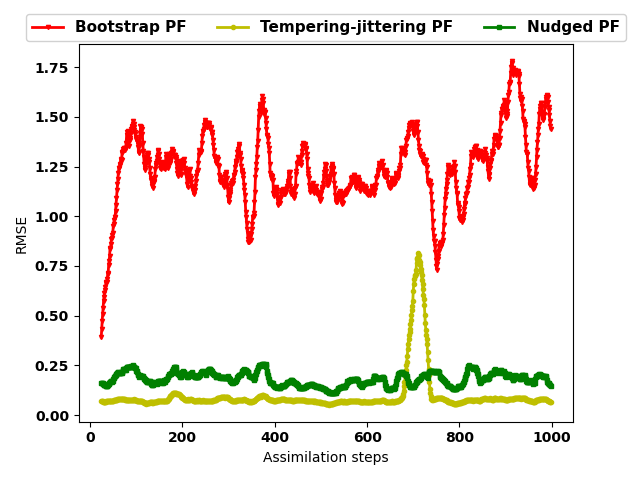} \\
        \includegraphics[width=0.5\linewidth]{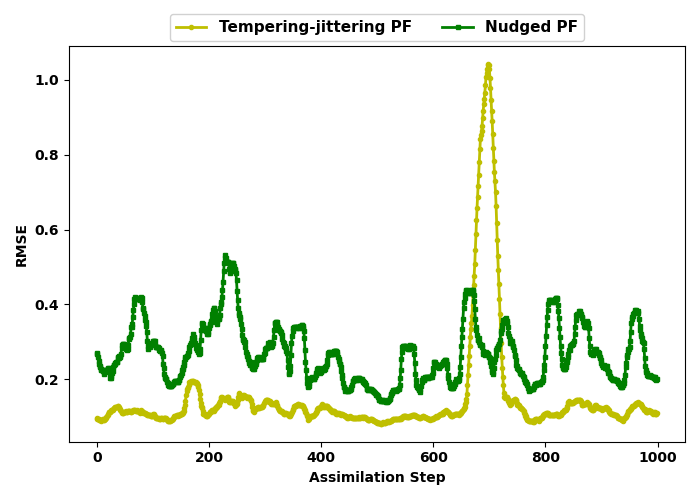}
        \caption{RMSE plots for the bootstrap, temper-jitter and Girsanov nudged filters applied to the SKS equation. Top: all three plotted. Bottom: y-axis rescaled for just the latter two filters.}
        \label{fig:RMSE}
\end{figure}
\begin{figure}[htb!]
        \includegraphics[width=.5\linewidth]{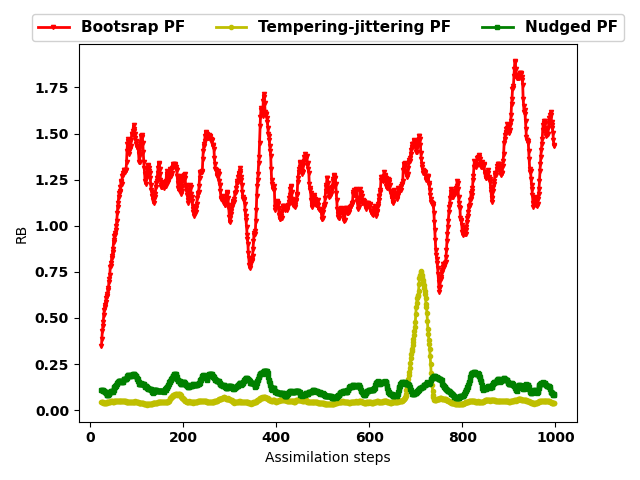} \\
        \includegraphics[width=.5\linewidth]{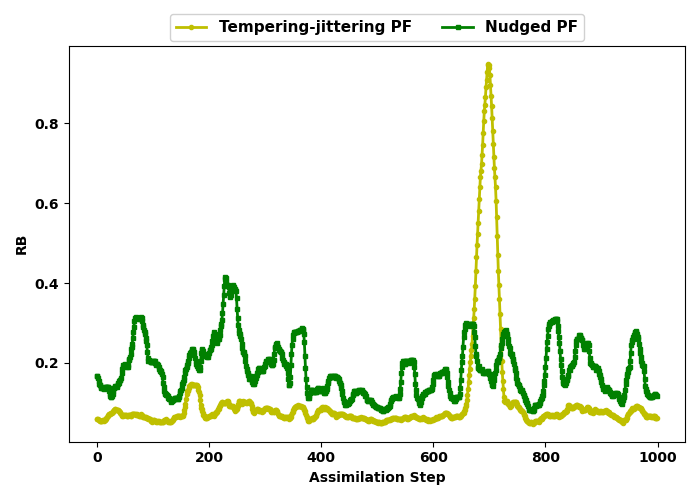} 
        \caption{RB plots for the bootstrap, temper-jitter and Girsanov nudged filters applied to the SKS equation.
        Top: all three plotted. Bottom: y-axis rescaled for just the latter two filters.}
        \label{fig:RB}
\end{figure}
\begin{figure}[htb!]
        \centering
        \includegraphics[width=.5\linewidth]{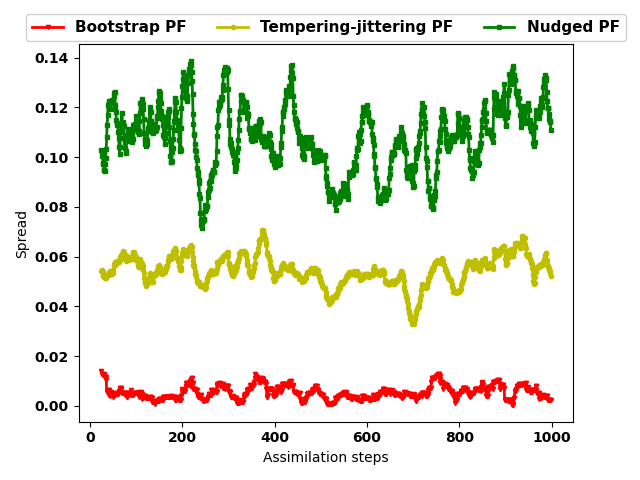} \\
        \includegraphics[width=.5\linewidth]{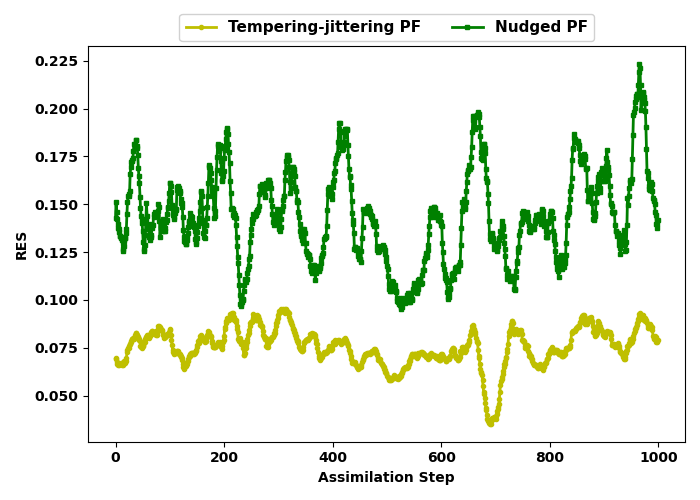} 
        \caption{(RES plots for the bootstrap, temper-jitter and Girsanov nudged filters applied to the SKS equation.
        Top: all three plotted. Bottom: y-axis rescaled for just the latter two filters.}
            \label{fig:RES}
            \end{figure}

\begin{figure}[htb!]
    \centering
    \begin{minipage}{0.48\textwidth}
        \centering
        \includegraphics[width=0.9\linewidth]{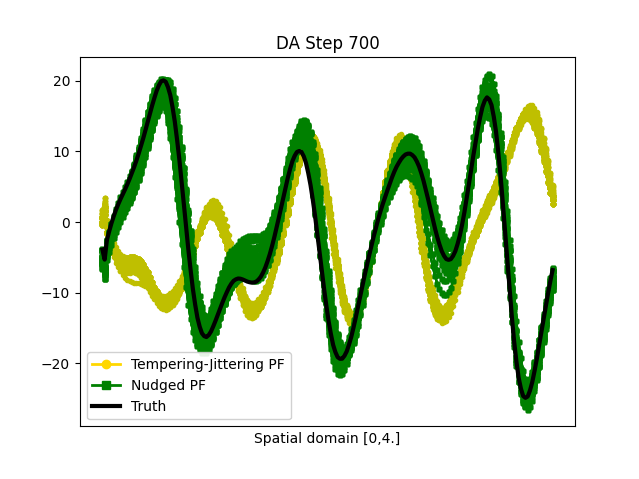}
        \subcaption{DA Step 700}
    \end{minipage}
    \hfill
    \begin{minipage}{0.48\textwidth}
        \centering
        \includegraphics[width=0.9\linewidth]{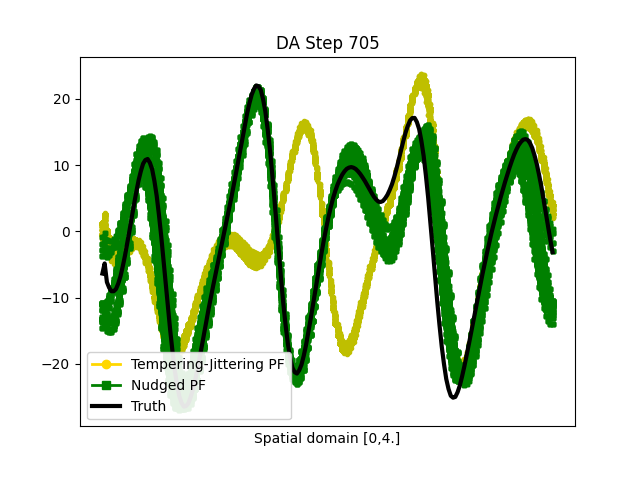}
        \subcaption{DA Step 705}
    \end{minipage}
    \vspace{0.5em}  
    \begin{minipage}{0.48\textwidth}
        \centering
        \includegraphics[width=0.9\linewidth]{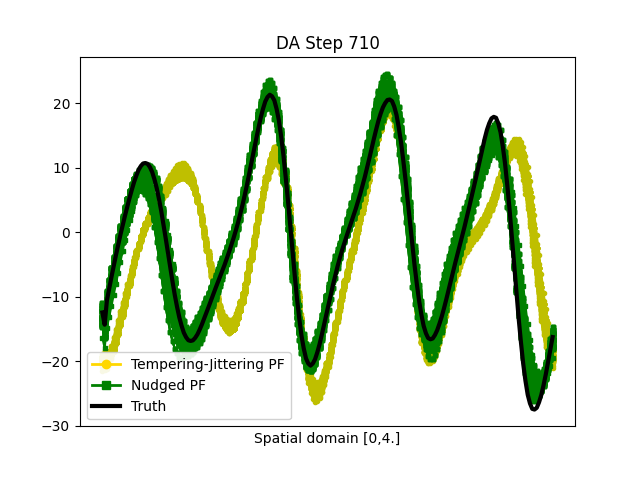}
        \subcaption{DA Step 710}
    \end{minipage}
    \hfill
    \begin{minipage}{0.48\textwidth}
        \centering
        \includegraphics[width=0.9\linewidth]{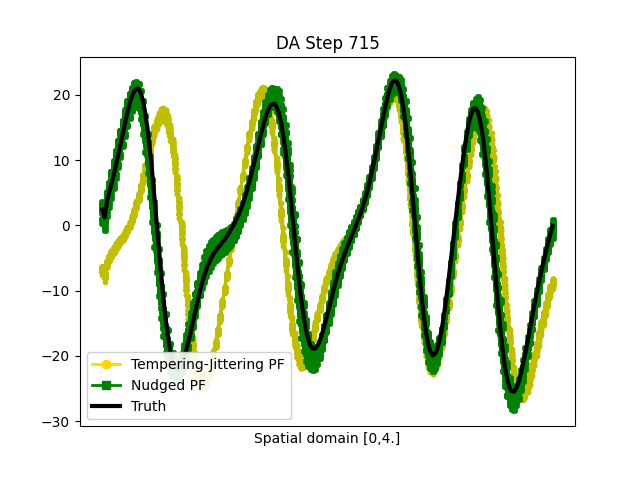}
        \subcaption{DA Step 715}
    \end{minipage}

    \caption{Comparison of ensemble vs truth for temper-jitter PF vs Nudged PF for various d.a. step}
    \label{fig:fdallerror_comparison_2x2}
\end{figure}

\begin{figure}[htb!]
    \centering
    \includegraphics[width=0.6\linewidth]{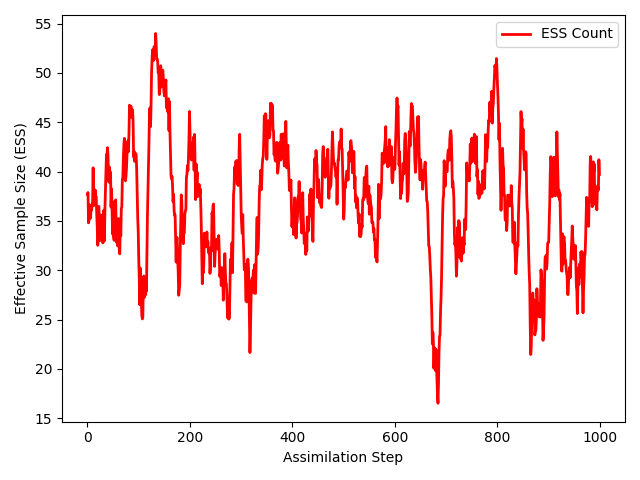}
    \caption{ESS}
    \label{fig:ess_plot}
\end{figure}

To quantify the behaviour of the particle filters, we consider the following widely used metrics,
\begin{align}
\mbox{RMSE} & := \dfrac{1}{N_p} \sum_{i=1}^{N_p}\dfrac{\|y_{true} - h(X_i)\|_{\ell^2}}{\|y_{true}\|_{\ell^2}}, \\
\mbox{RB} &:=  \dfrac{\|y_{true} - \bar{h(X)}\|_{\ell^1}}{\|y_{true}\|_{\ell^1}}, \qquad \mbox{where }
\bar{h(X)} = \frac{1}{N_p}\sum_{i=1}^{N_p}h(X_i),
\\
\mbox{RES} & :=  \dfrac{1}{N_p - 1} \sum_{p=1}^{N_p}\dfrac{\| \bar{h(X)}-h(X_i)\|^2_{\ell^2}}{\|y_{true}\|^2_{\ell^2}}.
\end{align}
Time series for these metrics are shown in Figures \ref{fig:RMSE}, \ref{fig:RB} and \ref{fig:RES}.
We observe that both the temper-jitter and Girsanov nudged particle filters are able to track the true solution in 
a stable manner. whilst the bootstrap filter is not (as expected). We also observe that the Girsanov nudged particle filter is producing an ensemble which is more spread than the temper-jitter filter. This appears to be because the former filter moves the particles in regions of highest likelihood to regions of lower likelihood to try to balance out the weights of the ensemble, whilst the latter filter repeatedly resamples during the tempering procedure to produce an ensemble where most of the particles originate from the same state at the previous assimilation step.
However, we also observe that the temper-jitter filter undergoes an excursion at around 700 data assimilation steps, where the ensemble is drifting away from the truth; at later times this excursion is arrested and the relative bias decreases again. This appears to occur because of an ``extreme event'' in the true dynamics, leading to a solution that
is unlikely with respect to the filtering distribution before this point, perhaps due to the temporary creation or destruction of a bump in the solution. We do not observe this large excursion with the Girsanov nudged filter, suggesting that it was able to move the ensemble closer to the the region of high likelihood using the control variables. We see some evidence for this in Figure \ref{fig:fdallerror_comparison_2x2}, which shows that the true solution undergoes a merger of two bumps into one, which does not occur in any of the temper-jitter states. The Girsanov nudged filter ensemble particles all undergo this merge after a shorter time when the true solution is outside the spread of the ensemble in that region of the domain. Figure \ref{fig:ess_plot} shows the time series
for ESS for the Girsanov nudged filter during the experiment. We see that it fluctuates a lot but it remains a significant fraction of the total ensemble size 90 throughout, with the lowest value during the extreme event studied in Figure \ref{fig:fdallerror_comparison_2x2}, as might be expected.

As a final comparison, we examine rank histograms for the particle filters. These are computed following the procedure described in \emph{e.g.} \cite{reich2015probabilistic} (Section 4.4). As expected, the bootstrap
filter produces extremely poor results and is not shown. We observe that both the temper-jitter and Girsanov nudging
filters behave quite similarly, with flat interiors but significant spikes at the sides of the histogram which are 
probably due to the excursions where the true value moves outside the ensemble.

\begin{figure}[htb!]
    \centering
        \includegraphics[width=0.65\linewidth]{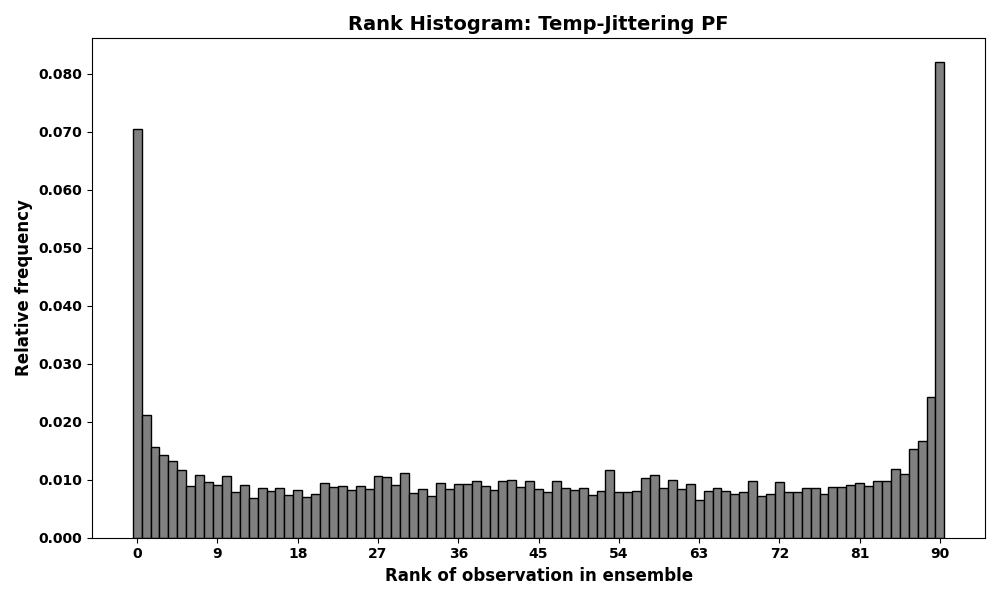} \\
        \includegraphics[width=0.65\linewidth]{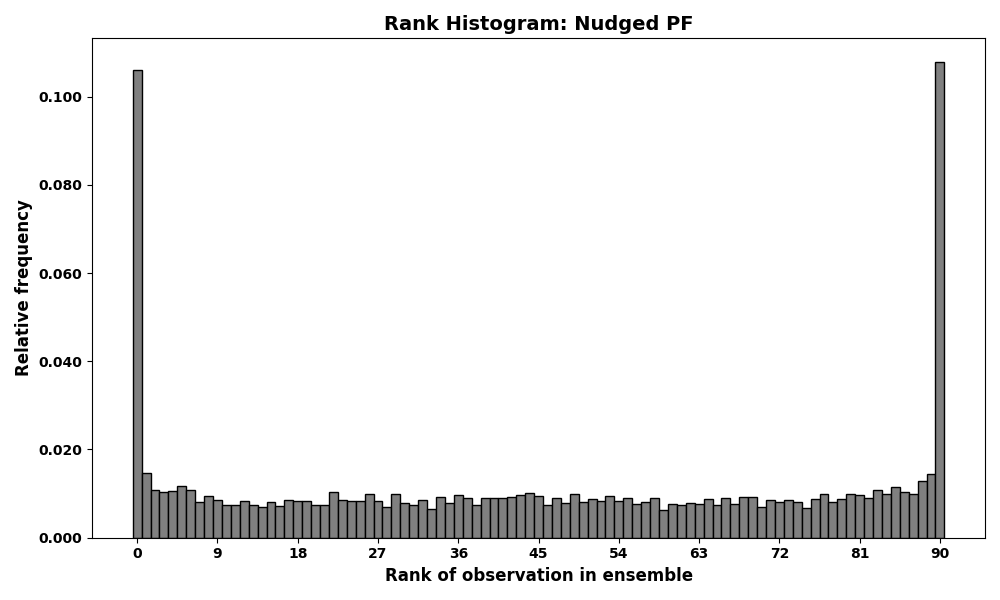} 
    \caption{Comparison of rank histograms. Top: the temper-jitter filter. Bottom: the Girsanov nudging filter.}
    \label{fig:rankHist}
\end{figure}

\section{Summary and outlook}
\label{sec:summary}

We presented a nudging approach to particle filters for SPDEs, based on introducing a Girsanov perturbed model. This allows the choice of control variables which can be used to steer the ensemble towards regions of higher likelihood. We propose that the control variables need to be optimised collectively in order to avoid low ESS and 
consequent filter divergence. We introduced a formulation of an optimisation problem that allows the control variables to be solved in three stages; these three stages separate the nonlinearity of the ESS formula and the nonlinearity of the forward model, and allow computations in the control variable space to be made in parallel
across the ensemble of particles. We presented numerical results that demonstrate the behaviour of the filter,
using our parallelised implementation developed using Firedrake, making specific use of its automated adjoint
capabilities for gradient based optimisation. In our comparisons, we used the temper-jitter filter as a scalable
reference, since it has been rigorously analysed and tested on numerous SDE and SPDEs. We observed that for the 
same ensemble size of 90, the Girsanov nudging filter showed a larger ensemble spread than the temper-jitter filter, with a larger relative bias in general. However, the Girsanov filter showed signs of additional stability, recovering
much more quickly after an extreme event when the true solution became unlikely with respect to the filtering distribution. This demonstrates the potential of Girsanov nudging for building more robust particle filters for PDEs.

It is clear that further work is needed to develop the filter, so that higher values of ESS can be achieved. We plan to develop extensions of the filter that interleave nudging with tempering and jittering, in the hope of combining the stability of one filter with the accuracy of the other. More engineering work is also needed to tune optimisation parameters in order to avoid wasting iterations (each of which involves solving the forward model) when they do not contribute to the stability or accuracy of the filter. We also plan to apply these particle filters to larger scale 
SPDEs using high performance computing facilities, bringing our work closer to applications.

\paragraph{\bfseries Acknowledgements} The authors acknowledge funding from EPSRC grant EP/W016125/1 "Next generation particle filters for stochastic partial differential equations", and the European Research Council Synergy grant "Stochastic Transport in Upper Ocean Dynamics".

\bibliographystyle{plain} 
\bibliography{sch_da}

\end{document}